\documentclass[12pt]{amsart}

\usepackage{amsmath,amssymb,amsthm,wrapfig,amsfonts}

\newtheorem{defn}{Definition}[section]
\newtheorem{thm}[defn]{Theorem}

\newtheorem{cor}[defn]{Corollary}
\newtheorem{lem}[defn]{Lemma}

\newtheorem{ack}{Acknowledgements\!\!}

\numberwithin{equation}{section}

\title[Approximated Brunn--Minkowski inequalities]
{Concentration of measure via approximated Brunn--Minkowski inequalities}

\author[Masayoshi Watanabe]{Masayoshi Watanabe}
\address{Mathematical Institute, Tohoku University, Sendai, 980-8578, Japan}
\email{sa3m33@math.tohoku.ac.jp}
\date{\today}

\begin{document}

\maketitle

\begin{abstract}
We prove that an approximated version of
the Brunn--Minkowski inequality with volume distortion coefficient
implies a Gaussian concentration-of-measure phenomenon.
Our main theorem is applicable to discrete spaces.
\end{abstract}

\section{Introduction}

\noindent
Let $(X,d)$ be a complete separable metric space 
equipped with
a Borel probability measure $\mu$ on $X$ with full support.
Henceforth, we call such a triple a \textit{metric measure space}.
The \textit{concentration function} of a metric measure space $(X,d,\mu)$ is defined by 
$$\alpha_{(X,\,d,\,\mu)}(r)=\sup \{ 1-\mu(A_r)\,|\,
A\mbox{ is a Borel set in }X\mbox{ with }\mu(A)\ge 1/2\},$$
where 
$A_r:=\{x\in X\,|\,d(x,A)<r\}$.

Let $(M,g)$ be an $n$-dimensional complete Riemannian manifold
with Riemannian distance $d_g$ and normalized Riemannian measure $\mu_g$.
If the Ricci curvature of $M$ is bounded below by $n-1$,
then the L\'{e}vy--Gromov isoperimetric inequality~\cite[Appendix C]{G} implies
\begin{equation*}\label{conc}
\alpha_{(M,\,d_g,\, \mu_g)}(r)\le e^{-(n-1)r^2/2}
\end{equation*}
for every $r>0$.
This is an example of \textit{Gaussian concentration-of-measure phenomenon}.
See \cite{Le} for details.
Moreover, the curvature-dimension condition CD$(n-1,n)$, or $n$-Ricci curvature $\ge n-1$,
for measured length spaces
implies a Gaussian concentration
via log-Sobolev inequality; 
see \cite[Corollary 6.12]{LVricci}, \cite{LVweak}, and \cite[Theorem 5.3]{Le}.

In this paper
we deduce
a Gaussian concentration
from
a weaker condition:
an $\epsilon$-\textit{approximated Brunn--Minkowski inequality} 
$\epsilon$-$\mathrm{BM}(n-1,n)$ \textit{of dimension} $n$ \textit{and of Ricci curvature} $\ge n-1$,
introduced by Bonnefont~\cite{B}.
The definition of $\epsilon$-$\mathrm{BM}(n-1,n)$ is in Section 2.
Our main theorem is

\begin{thm}\label{main}
Let $\epsilon \ge 0$ and $n\in (1,\infty)$.
If a metric measure space $(X,d,\mu)$ satisfies $\epsilon$-$\mathrm{BM}(n-1,n)$,
then we have, for every $r>0$,
$$\alpha_{(X,\,d,\,\mu)}(r)\le 2e^{-(n-1)r^2/\pi^2}.$$
\end{thm}

Note that the curvature-dimension condition CD$(n-1,n)$ does not make sense in discrete spaces;
however, $\epsilon$-$\mathrm{BM}(n-1,n)$ does.
For example, we can apply Theorem~\ref{main} to the discretization of a measured length space:
let us explain what we mean by the discretization in Section 4.

\begin{cor}\label{discretized}
Given $\epsilon \ge 0$ and $n\in (1,\infty)$,
let $(X_{\epsilon},d,\mu_{\epsilon})$ be a discretization
of a measured length space with 
the curvature-dimension condition $\mathrm{CD}(n-1,n)$.
Then we have, for every $r>0$,
$$\alpha_{(X_{\epsilon},\,d,\,\mu_{\epsilon})}(r)\le 2e^{-(n-1)r^2/\pi^2}.$$
\end{cor}

\begin{ack}\upshape
The author is sincerely grateful to Professor Takashi Shioya 
for his advice and constant encouragement.
\end{ack}

\section{Approximated Brunn--Minkowski inequality}

\noindent
Let $(X,d,\mu)$ be a metric measure space.
Given $\epsilon \ge 0$, $t\in (0,1)$, and $A_0,A_1\subset X$,
we first define
the set of $\epsilon$-\textit{approximated} $t$-\textit{intermediate points} between $A_0$ and $A_1$ 
by
\begin{align*}
I_t^{\epsilon}(A_0,A_1)=\big\{ &x\in X \bigm| \mbox{there exist } 
                               x_0\in A_0\mbox{ and } x_1\in A_1 \mbox{ with}\\
                               &|d(x_0,x)-td(x_0,x_1)|\le \epsilon \mbox{ and }
                                              |d(x,x_1)-(1-t)d(x_0,x_1)|\le \epsilon  \big\}.
\end{align*}                                                
In Euclidean space,
$I_t^0(A_0,A_1)$ coincides with the Minkowski sum $(1-t)A_0+tA_1$.

\begin{defn}\label{distortedBM}\upshape
Let $\epsilon \ge 0$ and $n\in (1,\infty)$.
We say that $(X,d,\mu)$ satisfies 
an $\epsilon$-\textit{approximated Brunn--Minkowski inequality of dimension} $n$ 
\textit{and of Ricci curvature} $\ge n-1$ or,
for short, $\epsilon$-$\mathrm{BM}(n-1,n)$
if we have
\begin{align}\label{BM}
\mu(I_t^{\epsilon}(A_0,A_1))^{1/n}\ge
(1-t)&\Bigg[\inf_{x_0\in A_0,\,x_1\in A_1}
\bigg(\frac{\sin((1-t)d(x_0,x_1))}{(1-t)\sin d(x_0,x_1)}\bigg)^{(n-1)/n}\Bigg]\mu(A_0)^{1/n}\\
&+t\Bigg[\inf_{x_0\in A_0,\,x_1\in A_1}
\bigg(\frac{\sin(td(x_0,x_1))}{t\sin d(x_0,x_1)}\bigg)^{(n-1)/n}\Bigg]\mu(A_1)^{1/n}\notag
\end{align}
for all nonempty Borel sets $A_0,A_1\subset X$
and for all $t\in (0,1)$,
where
\begin{equation*}
\frac{\sin(td(x_0,x_1))}{t\sin d(x_0,x_1)}:=+\infty \mbox{\ \ \ \ \ \rm{if}\ \ } d(x_0,x_1)\ge \pi.
\end{equation*}
\end{defn}

See \cite[Section 14]{Voptimal} for the meaning of distortion coefficients in (\ref{BM}).
Clearly, $\epsilon$-BM$(n-1,n)$ implies $\epsilon^{\prime}$-BM$(n-1,n)$
for $\epsilon^{\prime}\ge \epsilon$.
The curvature-dimension condition CD$(n-1,n)$ 
implies $0$-$\mathrm{BM}(n^{\prime}-1,n^{\prime})$ for all $n^{\prime}\ge n$;
see \cite[Proposition 2.1]{SII} and \cite[Theorem 30.7]{Voptimal}.
The Brunn--Minkowski inequality in curved spaces
is proved by virtue of \cite{CMS}.

\section{Concentration of measure}

\noindent
We begin with a lemma
corresponding to the Bonnet--Myers theorem;
see \cite[Corollary 2.6]{SII} and \cite[Proposition 29.11]{Voptimal}.

\begin{lem}\label{diam}
Let $\epsilon \ge 0$ and $n\in (1,\infty)$.
If a metric measure space $(X,d,\mu)$ satisfies $\epsilon$-$\mathrm{BM}(n-1,n)$,
then
$\mathrm{diam} (X)\le \pi$.
\end{lem}
\begin{proof}
Suppose that there are two points $x_0,x_1\in X$
with $d(x_0,x_1)>\pi$.
Choosing a sufficiently small $\delta>0$,
we have $d(B_{\delta}(x_0),B_{\delta}(x_1))>\pi$.
Note that, in (\ref{BM}) with
$A_0=B_{\delta}(x_0)$ and $A_1=B_{\delta}(x_1)$,
the coefficients in the right-hand side equals $+\infty$.
We  then have a contradiction from
$\mu(A_0)>0$, $\mu(A_1)>0$,
and $\mu(I_t^{\epsilon}(A_0,A_1))\le \mu(X)=1$.
\end{proof}

\begin{proof}[Proof of Theorem~\ref{main}]
Let $A$ be a Borel set in $X$ with $\mu(A)\ge 1/2$.
By Lemma~\ref{diam},
it suffices to prove
$1-\mu(A_r)\le 2e^{-(n-1)r^2/\pi^2}$
for every $r\in (0,\pi)$.

We now put $B=X\setminus A_r$
for a fixed $r\in (0,\pi)$.
Note that $(\sin(d/2))/((1/2)\sin d)$ is
monotone nonincreasing in $d\in (0,\pi)$.
Since $d(A,B)\ge r$,
it follows that
$$\inf_{x\in A,\,y\in B}
\bigg(\frac{\sin((1/2)d(x,y))}{(1/2)\sin d(x,y)}\bigg)^{(n-1)/n}
\ge \bigg(\frac{\sin(r/2)}{(1/2)\sin r}\bigg)^{(n-1)/n}.$$
Inequality (\ref{BM}) with $A_0=A$, $A_1=B$, and $t=1/2$ gives
\begin{align*}
\mu(Z_{1/2}^{\epsilon}(A,B))^{1/n}
&\ge 
\frac{1}{2}\bigg(\frac{\sin (r/2)}{(1/2)\sin r}\bigg)^{(n-1)/n}
(\mu(A)^{1/n}+\mu(B)^{1/n})\\
&\ge \frac{(\mu(A)^{1/n}\mu(B)^{1/n})^{1/2}}{(\cos (r/2))^{(n-1)/n}}.
\end{align*}
We used relations
$\sin r=2\sin(r/2)\cos(r/2)$ and
$(a+b)/2\ge \sqrt{ab}$
for the last step.
Noting 
$\mu(A)\ge 1/2$ and
$\mu(Z_{1/2}^{\epsilon}(A,B))\le \mu(X)=1$,
we get
$$\mu(B)\le 2\Big(\cos \frac{r}{2}\Big)^{2(n-1)}
\le 2\Big(1-\frac{r^2}{2\pi^2}\Big)^{2(n-1)}
\le 2e^{-(n-1)r^2/\pi^2}.$$
\end{proof}


We can get a better estimate 
for all sufficiently large $n \in (1,\infty)$ and
small $r>0$.

\begin{thm}
Fix $\epsilon \ge 0$.
Given $\delta>0$,
there exist $n_0 \in (1,\infty)$ and $r_0>0$
such that
if a metric measure space $(X,d,\mu)$ satisfies $\epsilon$-$\mathrm{BM}(n-1,n)$
for a number $n\ge n_0$,
then we have,
for $0<r\le r_0$,
$$\alpha_{(X,\,d,\,\mu)}(r)\le e^{-(1-\delta)nr^2/4}.$$
\end{thm}

\begin{proof}
Modify the proof of Theorem~\ref{main}.
We get
$$\alpha_{(X,\,d,\,\mu)}(r)\le e^{-n[1+2^{-1/n}-2(\cos(r/2))^{(n-1)/n}]}$$
without employing the arithmetic-geometric mean inequality.
Taylor expansion, $2-2\cos(r/2)=r^2/4+o(r^2)$, completes the proof.
\end{proof}

\section{Discretization}

\noindent
Let $(X,d,\mu)$ be a metric measure space.
Given $\epsilon>0$,
take a set $\{x_i\}_{i=1}^{\infty}$
of countable distinct points in $X$ with $X=\bigcup_{i=1}^{\infty}B_{\epsilon}(x_i)$,
where $B_{\epsilon}(x_i)$ is the open ball
of radius $\epsilon$ centered at $x_i$.
We can choose a measurable set $A_i\subset B_{\epsilon}(x_i)$
for each $i$ such that
$x_i\in A_i$, $A_i\cap A_j\neq \emptyset$ $(i\neq j)$, and $X=\bigcup_{i=1}^{\infty}A_i.$
Setting $\mu_{\epsilon}(\{ x_i\})=\mu(A_i)$,
we get a probability measure $\mu_{\epsilon}$ on $X_{\epsilon}:=\{x_i\}_{i=1}^{\infty}$.
We call $(X_{\epsilon},d,\mu_{\epsilon})$ 
a \textit{discretization} of $(X,d,\mu)$.

\begin{proof}[Proof of Corollary~\ref{discretized}]
Every discretization $(X_{\epsilon},d,\mu_{\epsilon})$
of a measured length space $(X,d,\mu)$ with 
CD$(n-1,n)$ satisfies
$4\epsilon$-BM$(n-1,n)$ \cite[Section 3]{B};
therefore,
Theorem~\ref{main} completes the proof.
\end{proof}



\begin{thebibliography}{99}
\bibitem{B}
 M. Bonnefont,
 A discrete version and stability of Brunn Minkowski inequality,
 preprint (2007).
\bibitem{CMS}
 D. Cordero-Erausquin, R. J. McCann, and M. Schmuckenschl\"{a}ger,
 A Riemannian interpolation inequality \`{a} la Borell, Brascamp and Lieb,
 Invent. Math. {\bf 146} (2001), no. 2, 219--257. 
\bibitem{G}
 M. Gromov,
 \textit{Metric structures for Riemannian and non-Riemannian spaces},
 Modern Birkh\"{a}user Classics, Birkh\"{a}user Boston, Inc., Boston, MA (2007). 
\bibitem{Le} 
 M. Ledoux,
 \textit{The concentration of measure phenomenon},
 Mathematical Surveys and Monographs {\bf 89},
 American Mathematical Society, Providence, RI (2001).
\bibitem{LVricci}
 J. Lott and C. Villani,
 Ricci curvature for metric-measure spaces via optimal transport, 
 Ann. of Math. (to appear).
\bibitem{LVweak}
 J. Lott and C. Villani,
 Weak curvature conditions and functional inequalities,
 J. Funct. Anal. {\bf 245} (2007), no. 1, 311--333.
\bibitem{SI}
 K.-T. Sturm,
 On the geometry of metric measure spaces I, 
 Acta Math. {\bf 196} (2006), no. 1, 65--131. 
\bibitem{SII}
 K.-T. Sturm,
 On the geometry of metric measure spaces II, 
 Acta Math. {\bf 196} (2006), no. 1, 133--177.
\bibitem{Voptimal}
 C. Villani,
 \textit{Optimal transport, old and new},
 Grundlehren der Mathematischen Wissenschaften, Springer-Verlag
 (to appear).
\end{thebibliography}
\end{document}